\documentstyle[12pt]{article} 
\newcommand{\doublespace}{
   \renewcommand{\baselinestretch}{1.2}
   \large\normalsize}

\def \Z{\Bbb Z}
\def \C{\Bbb C}

\def \End{{\rm End}}
\def \Hom{{\rm Hom}}

\def \<{\langle} 
\def \>{\rangle} 

\def \a{\alpha }

\def \l{\lambda }
\def \L{\Lambda }

\def \b{\beta }

\def \s{\sigma}

\def \qed{\mbox{ $\square$}}
\def \pf {\noindent {\bf Proof:} \,}
\input amssym.def
\input amssym
\doublespace
\begin{document}
\newtheorem{thmm}{Theorem}
\newtheorem{thm}{Theorem}[section]
\newtheorem{prop}[thm]{Proposition}
\newtheorem{cor}[thm]{Corollary}
\newtheorem{lem}[thm]{Lemma}
\newtheorem{rem}[thm]{Remark}
\newtheorem{de}[thmm]{Definition}
\begin{center}
{\Large {\bf Quantum Galois theory for compact Lie groups}} \\
\vspace{0.5cm}
Chongying Dong\footnote{Supported by NSF grant 
DMS-9700923 and a research grant from the Committee on Research, UC Santa Cruz.} and Geoffrey Mason\footnote{Supported by NSF grant DMS-9700909 and a research grant from the Committee on Research, UC Santa Cruz.}
\\
Department of Mathematics, University of
California, Santa Cruz, CA 95064
\end{center}

\hspace{1.5 cm}
\begin{abstract} We establish a quantum Galois correspondence for compact
Lie groups of automorphisms acting on a simple vertex operator algebra.
\end{abstract}

\section{Introduction}

Suppose that $V$ is a simple vertex operator algebra (VOA) and that
$G$ is a compact Lie group, possibly finite, which acts as a continuous 
group of automorphisms of $V.$ Such pairs $(V,G)$ occurs in a number of 
contexts within conformal field theory, for example the theory
of orbifolds [DM1], and the theory of $W$-algebras [BS], and
it is of interest to investigate the general situation.

In [DM2] we introduced some techniques for studying the representation
of $V^G,$ the sub VOA of fixed points of $G$ on $V,$ afforded by $V,$ 
and used this to establish a quantum Galois correspondence 
for certain classes of finite groups $G.$ Later, in [DLM], the 
representation of $V^G$ on $V$ was considered for a general compact Lie group $G$ and the quantum Galois correspondence was 
established for compact {\em abelian} groups. Very recently, Hanaki, Miyamoto and Tambara have considered the correspondence for an {\em arbitrary} 
finite group [HMT].

The purpose of present paper is to establish a quantum Galois
correspondence for a general compact Lie group. To describe the main
result we need to recall the following result (Theorem 2.4 of [DLM])
which is fundamental to our approach.
\begin{thmm}\label{t2.1} There is a decomposition of $V,$ considered as
$G\times V^G$-module,
\begin{equation}\label{2.3}
V=\bigoplus_{\l\in\L} W_{\l}\otimes V_{\l}.
\end{equation}
Here, $\Lambda$ is  a set which indexes the inequivalent, finite-dimensional,
continuous, complex, simple, unitary representations of $G,$ and
$W_{\l}$ is the corresponding left $G$-module. In the decomposition
(\ref{2.3}), the following hold:

(a) For each $\l\in\L,$ $V_{\l}=\Hom_G(W_{\l},V)$ is a {\em nonzero, irreducible} $ V^G$-module.

(b) If $\l, \mu\in \L$ and $\l\ne \mu$ then $V_{\l}$ and $V_{\mu}$ are 
{\em not} isomorphic as $V^G$-modules.
\end{thmm}

Now let $U$ be sub VOA of $V$  which contains $V^G.$ From (\ref{2.3}) we see that
\begin{equation}\label{n1.2}
U=\bigoplus_{\l\in \L}R_{\l}\otimes V_{\l}
\end{equation}
for subspaces $R_{\l}\subset W_{\l}.$ Let $\<\cdot,\cdot\>_{\l}$ 
be a positive definite, $G$-invariant, hermitian form on $W_{\l},$
and let $R_{\l}^c$ be the orthogonal complement to $R_{\l}$ in
$W_{\l}.$ 

\begin{de} {\rm We say that $U$ is {\em orthogonally complemented}
in $V$ if the subspace $U^c$ defined by
\begin{equation}\label{n1.3}
  U^c=\bigoplus_{\l\in \L}R_{\l}^c\otimes V_{\l}
\end{equation}
is a $U$-submodule  of $V.$}
\end{de}

Our main result is the following:
\begin{thmm}\label{thmm} Let $V$ be a simple vertex operator algebra and let $G$ be a 
compact Lie group  of automorphisms of $V$ which acts continuously on $V.$ Let
\begin{equation}\label{1.1}
\gamma: H\mapsto V^H
\end{equation}
be the map which associates to a subgroup $H$ of $G$ the sub VOA $V^H$ consisting
of $H$-fixed points. Then the following hold:

(A) $\gamma$ induces a bijection between the {\em closed} Lie subgroups
of $G$ and the {\em simple,} {\em orthogonally complemented} sub VOAs of $V.$ 

(B) If $G$ is finite then $\gamma$ induces a bijection between the subgroups 
of $G$ and the sub VOAs of $V$ which contain $V^G.$
\end{thmm}

Part (A) is the main new result; (B) is a restatement of the result
of Hanaki, Miyamoto and Tambara. Of course, the classical Galois correspondence 
may be stated in the following form: if $K$ is a commutative 
field and $G$ a (finite) group of automorphisms of $K$ then 
$K$ is a Galois extension of $K^G$ and $H\mapsto K^H$ set-up a bijection
between subgroups of $G$ and subfields of $K$ which contain $K^G.$ 
In this sense the theorem represents a direct generalization of the classical 
case inasmuch as a commutative field may be considered to be a simple 
VOA. 

One of the new ingredients appearing in (A) is the fact that one has to restrict 
oneself to {\em simple} sub VOAs of $V$ which contain $V^G,$ a condition 
that is obviously redundant in the classical case. Notice that the
adjectives ``simple'' and ``orthogonally complemented'' do not appear in (B).  
Indeed, given 
the truth of (A), (B) is equivalent to the assertion that, if $G$ is finite,
{\em every} sub VOA of $V$ which contains $V^G$ is simple and orthogonally complemented.  
It should
be noted that the simplicity condition
appearing in (A) is necessary: there are many examples
of sub VOAs of $V$ which contain $V^G$ which are not simple. We discuss 
this and related matters in Section 4. Such sub VOAs cannot 
appear in a Galois correspondence such as ours. It is likely that 
the complementary condition in (A) is unnecessary, but we do not
know how to remove it at present.

Our proof follows in broad outline that of [DLM] with the new ingredient supplied 
by [HMT]  added in a suitably modified way. We use a basic result
from conformal field theory (established in Section 2) together
with some standard techniques from Lie group theory (representative functions
and the Stone-Weierstrass Theorem).

We thank Hanaki, Miyamoto and Tambara for communicating their pretty paper
[HMT] and thereby inspiring the results contained herein.

\section{Some results from VOA theory}
\setcounter{equation}{0}

The following notation will be in force throughout  of the paper:
$V=(V,Y,{\bf 1},\omega)$ is a simple vertex operator algebra and
$G\leq Aut(V)$ a compact Lie group of automorphisms of $V$ which acts
 continuously on $V.$ For a subgroup $H\leq G,$ $V^H$ is the sub VOA of
$H$-invariants and $\gamma$ the map (\ref{1.1}). We retain the other
notation introduced in Section 1.

Recall that $V$ is a direct sum of finite-dimensional homogeneous subspaces
\begin{equation}\label{2.1}
V=\bigoplus_{n\in \Z}V_n
\end{equation}
where $V_n=0$ for all small enough $n.$ If $v\in V$ the corresponding
vertex operator is denoted by
\begin{equation}\label{2.2}
Y(v,z)=\sum_{n\in\Z}v_nz^{-n-1}
\end{equation}
so that $v_n\in\End (V).$ (For basic facts about VOAs, we refer the
reader to [FHL].)

We may, and shall, choose each $W_{\l}$ to be homogeneous i.e., 
$W_{\l}\subset V_k$ for some $k\in \Z.$ This is possible because each $V_k$
is a $G$-module.

Now let $U$ be a sub VOA of $V$ which satisfies
\begin{equation}\label{2.4}
V^G\subset U\subset V,
\end{equation}
with decomposition of $U$ into $V^G$-modules as in (\ref{n1.2}) 
where $R_{\l}\subset W_{\l}.$ We set
\begin{equation}\label{2.6}
W=\bigoplus_{\l\in \L}W_{\l}
\end{equation}
\begin{equation}\label{2.7}
R=\bigoplus_{\l\in \L}R_{\l}
\end{equation}
and consider $W$ as a $G$-submodule of $V.$ Of course $R\subset W.$

The main idea of [HMT] is to establish the following
result:

\begin{prop}\label{p2.2}
Suppose that $\pi$ is a $G$-homomorphism from $W\otimes_{\C}W$ to $W.$ Then
$\pi(R\otimes R)\subset R.$
\end{prop}

Fix $\l,\mu\in \L$ and consider the following diagram of $G$-modules
and $G$-homomorphisms:
\begin{equation}\label{2.8}
\begin{array}{ccccccc} 
\cdots &\longrightarrow & Z_n &\stackrel{\phi_n}{\longrightarrow}&Z_{n+1}&\longrightarrow & \cdots\\
      &    &\Big\uparrow\vcenter{% 
\rlap{$\psi_n$}}& & \Big\uparrow\vcenter{% 
\rlap{$\psi_{n+1}$}}& &  \\
& & W_{\l}\otimes W_{\mu} & = &   W_{\l}\otimes W_{\mu} & & 
\end{array}
\end{equation}
The essential ingredients in (\ref{2.8}) are defined as follows:
for each integer $n,$
\begin{equation}\label{2.9}
Z_n=\<\sum_{m=n}^{\infty}u_mv|u\in W_{\l}, v\in W_{\mu}\>
\end{equation}
  \begin{equation}\label{2.10}
\psi_n: u\otimes v\mapsto \sum_{m=n}^{\infty}u_mv, u\in W_{\l}, v\in W_{\mu}
\end{equation}
 \begin{equation}\label{2.11}
\phi_n: \sum_{m=n}^{\infty}u_mv\mapsto \sum_{m=n+1}^{\infty}u_mv, u\in W_{\l}, v\in W_{\mu}.
\end{equation}
We make several observations. First, given $u$ and $v,$ we have $u_pv=0$
for all large enough $p.$ So all sums in (\ref{2.9})-(\ref{2.11})
are finite and $Z_n$ is a finite-dimensional $G$-submodule of $V.$  Second,
the maps $\psi_n$ and $\phi_n$ are well-defined. For $\psi_n$ this is clear;
as for $\phi_n,$ suppose that $\sum_{m=n}^{\infty}u_mv=0.$ Since $u$ and
$v$ are homogeneous (since $W_{\l},$ $W_{\mu}$ are so chosen), then $u_mv$ 
and $u_pv$ have different weights if $m\ne p.$ So $u_mv=0$ for all
$m\geq n$ and $\sum_{m=n+1}^{\infty}u_mv=0.$
Now the assertion concerning $\phi_n$ follows. Finally, note that 
indeed $\phi_n$ and $\psi_n$ are $G$-module homomorphisms. This
follows because if $g\in G$ then $gu_mg^{-1}=(gu)_m.$

\begin{lem}\label{l2.3}

For all small enough $n,$ $\psi_n$ is an isomorphism.
\end{lem}

\pf Assume false. Since each $\psi_n$ and $\phi_n$ is a surjection, there 
is an integer $q$ such that $\phi_n$ is an isomorphism for $n<q$ and moreover
there are $u^i\in W_{\l}, v^i\in W_{\mu},$ $i=0,...,t$ with
$v^0,...,v^t$ linearly independent and such that $$0\ne
\sum_{i=0}^tu^i\otimes v^i\in \ker \psi_q.$$ We thus have
$\sum_{i=0}^tu^i\otimes v^i\in \cap_{n\in\Z}\ker \psi_n,$ so for all
$n\in \Z$ we get $\sum_i\sum_{m=n}^{\infty}u^i_mv^i=0.$ But then for
each $n\in \Z$ we have $\sum_iu^i_nv^i=0,$ so that
\begin{equation}\label{2.12}
\sum_iY(u^i,z)v^i=0.
\end{equation}
Since the $v^i$ are linearly independent, (\ref{2.12}) forces each $u^i=0$
by Lemma 3.1 of [DM2], a contradiction because $\sum_iu^i\otimes v^i\ne 0.$
This completes the proof of the lemma. \qed

Turning to the proof of Proposition \ref{p2.2}, let $\pi$ be as in the statement of the proposition. It is sufficient to assume that $\pi$
induces a surjection $\pi: W_{\l}\otimes W_{\mu}\to W_{\nu}$ for
some $\l,\mu,\nu$ in $\L$ and establish that $\pi(R_{\l}\otimes R_{\mu})\subset R_{\nu}.$ There are $G$-module homomorphisms
$$
\begin{array}{ccccc}
Z_n&\stackrel{\psi_n^{-1}}{\longrightarrow} & W_{\l}\otimes W_{\mu}&
\stackrel{\pi}{\longrightarrow} & W_{\nu}\\
\Big\downarrow\vcenter{% 
\rlap{$i$}} & & &  & \Big\uparrow\vcenter{% 
\rlap{$\rho$}}\\
V & &= & & V
\end{array}
$$
where $i$ is inclusion, $n$ is chosen so that $\psi_n$ is an isomorphism (Lemma
\ref{l2.3}), and $\rho$ is some extension of $\pi\circ \psi_n^{-1}$ to 
$V.$

We claim that $\rho(U)\subset R.$ Indeed it is clear that $\rho$ annihilates 
$W_{\a}\otimes V_{\a}$ if $\a\ne \nu$ and induces a surjection 
$$\rho: W_{\nu}\otimes V_{\nu}\to W_{\nu}.$$
As $V_{\nu}$ is a trivial $G$-module, Schur's lemma tells us that if
$v\in V_{\nu}$ then the restriction of $\rho$ to
$W_{\nu}\otimes v$ has the form $\rho(w\otimes v)=kw$ for
$w\in W_{\nu}$ and $k$ a scalar. In particular $\rho(R_{\nu}\otimes V_{\nu})
\subset R_{\nu},$ so that $\rho(U)=\rho(R_{\nu}\otimes V_{\nu})\subset R$ as
claimed.

Finally, since $U$ is a sub VOA of $V$ we see that 
$\psi_n(R_{\l}\otimes R_{\mu})\subset U,$ so that
$$\pi( R_{\l}\otimes R_{\mu})=\rho\circ \psi_n(R_{\l}\otimes R_{\mu})
\subset \rho(U)\subset R.$$
This completes the proof of Proposition \ref{p2.2}.

\section{Proof of Theorem 1}
\setcounter{equation}{0}

We will use several standard results from the representation theory of
compact Lie groups. For these results and general background, the
reader may consult, for example, [BT].

It is a consequence of the Peter-Weyl theorem that if $H$ is a proper
closed subgroup of $G$ then there is $\l\in\L$ such that the space $W_{\l}^H$
of $H$-invariants is non-zero and $W_{\l}$ is not the trivial module.
 Since $W_{\l}$ is a constituent of $V$ 
by Theorem \ref{t2.1} it follows that $V^H\ne V^G.$ Therefore, the
map $\gamma$ (\ref{1.1}) is certainly an injection on the family of closed 
subgroups of $G.$

\begin{lem}\label{nl3.1} Let $H\subset G$ be a closed subgroup. Then $V^H$
is a simpleVOA which is orthogonally complemented in $V.$
\end{lem}

\pf The simplicity of $V^H$ follows from assertion (a) of Theorem \ref{t2.1}
applied with $H$ in place of $G;$ indeed, $V^H=V_{1_H}$ where $1_H$ is
the trivial $H$-module.

Now take $U=V^H$ as in (\ref{n1.2}), with $U^c$ as in (\ref{n1.3}). It
is easy to see that $R_{\l}^c$ is, in this case, the unique
$H$-invariant complement to $R_{\l}$ in $W_{\l},$ since $R_{\l}=W_{\l}^H.$
Thus, applying (\ref{2.3}) with $H$ in place of $G,$ we see that
\begin{equation}\label{n3.1}
U^c=\bigoplus_{\alpha}X_{\alpha}\otimes V_{\alpha}
\end{equation}
where $X_{\alpha}$ ranges over the non-trivial, simple, unitary $H$-modules,
and $V_{\a}$ is a simple $U$-module. In particular, $U^c$ is a 
$U$-module, so that $U$ is indeed orthogonally complemented. The lemma
is thus proved. \qed

It remains to prove that if $U$ is an orthogonally complemented simple
sub VOA of $V,$ then $U=V^H$ for some closed subgroup $H$ of $G.$ In the
course of proving this we will also establish (B) of Theorem
\ref{t2.1}.

We begin with an {\em arbitrary} sub VOA $U$ and assume earlier notation.
 Let $F$ be the space of complex-valued
representative functions on $G;$ $F$ is a $G$-bimodule via
$(g\a h)(k)=\a(g^{-1}kh^{-1})$ for $\a\in F,$ $g,h,k\in G,$ and there
is an isomorphism of $G$-bimodules
\begin{equation}\label{3.1}
F\cong\bigoplus_{\l\in \L}\End_{\C}(W_{\l})
\end{equation}
where $\End_{\C}(W_{\l})=W_{\l^*}\otimes W_{\l}$ carries the usual $G$-bimodule
structure.

Fix a basis for each $W_{\l}.$ Then  $\End_{\C}(W_{\l})$ may be identified
with a matrix algebra. Embed $W_{\l}$ into the first column of
$\End_{\C}(W_{\l^*}).$ This gives an embedding of $W$ into 
$F$ as a left $G$-submodule. Then $R$ is also a subspace of $F,$ and we let $S$ be the 
right $G$-submodule of $F$ generated by $R.$

Now $F$ is an algebra with respect to pointwise multiplication of
functions.  We consider this multiplication to be a map of
$G$-bimodules $\pi: F\otimes F\to F.$ We will establish
\begin{lem} $S$ is a subalgebra of $F.$
\end{lem}

\pf We must show that $\pi(S\otimes S)\subset S,$ and for this
it suffices to establish that $\pi(R_{\l}\otimes R_{\mu})\subset S$
for all $\l,\mu\in \L.$ Consider the diagram of left $G$-modules and maps
$$\begin{array}{ccc}
F\otimes F&\stackrel{\pi}{\longrightarrow}& F\\
 \Big\uparrow & & 
\Big\uparrow\\
\End\,W_{\l}\otimes\End\,W_{\mu} &\stackrel{\pi}{\longrightarrow}&\bigoplus_{\nu}\End\,W_{\nu}\\
 \Big\uparrow & & 
\Big\uparrow\\
W_{\l}\otimes W_{\mu} &\stackrel{\pi}{\longrightarrow}&\pi(W_{\l}\otimes W_{\mu})\\
& & \Big\downarrow\vcenter{% 
\rlap{$\rho'$}}\\
\Big\downarrow\vcenter{% 
\rlap{$\pi$}}
& & W_{\nu}'\\
& & \Big\downarrow\vcenter{% 
\rlap{$\cong$}}\\
\pi(W_{\l}\otimes W_{\mu}) &\stackrel{\rho}{\longrightarrow}& W_{\nu}
\end{array}$$
Here, $\nu$ ranges over those elements of $\L$ such that $W_{\nu}$ is a 
constituent of $W_{\l}\otimes W_{\mu},$ $\rho$ and $\rho'$ 
are projections onto $W_{\nu}$ and 
a summand $W_{\nu}'$ isomorphic to $W_{\nu},$ and the upper
maps are containments. By Proposition \ref{p2.2}
we obtain the following by restricting the bottom part of the preceding 
diagram:
$$\begin{array}{ccc}
 $$R_{\l}\otimes R_{\mu} &\stackrel{\pi}{\longrightarrow}&\pi(R_{\l}\otimes R_{\mu})\\
& & \Big\downarrow\vcenter{% 
\rlap{$\rho'$}}\\
\Big\downarrow\vcenter{% 
\rlap{$\pi$}}
& & R_{\nu}'\\
& & \Big\downarrow\vcenter{% 
\rlap{$\cong$}}\\
\pi(R_{\l}\otimes R_{\mu}) &\stackrel{\rho}{\longrightarrow}& R_{\nu}
\end{array}$$
Now the space of all such $W_{\nu}',$ that is to say $\End (W_{\nu}),$
is spanned by $W_{\nu}\cdot g$ for $g\in G,$ so the same is true
if we replace $W_{\nu}'$ by $R_{\nu}'$ and $W_{\nu}$ by $R_{\nu}.$ As
a result, we get
$$\pi(R_{\l}\otimes R_{\mu})\subset \bigoplus_{\nu}\sum_{g\in G}R_{\nu}\cdot g
\subset S.$$
The lemma is proved. \qed

We define
\begin{equation}\label{l3.2}
H=\{h\in G|\sigma (h)=\sigma(1), \forall \sigma\in S\}
\end{equation}
where $1$ refers to the identity of $G.$

Now obviously $1\in H.$ Remembering that $S$ is a right $G$-module,
if $h,k\in H$ and $\sigma\in S$ then
$$\sigma(hk)=(\sigma\cdot k^{-1})(h)=(\sigma\cdot k^{-1})(1)=\sigma(k)=\sigma(1)$$
so also $hk\in H.$ Similarly
$$\sigma(1)=\sigma(hh^{-1})=(\sigma\cdot h)(h)=(\sigma\cdot h)(1)
=\sigma( h^{-1})$$
so $h^{-1}\in H.$ This shows that $H$ is a closed subgroup of $G,$ in particular $H\backslash G$ is a Hausdorff space.

Next, if $\s\in S,$ $h\in H$ and $t\in G$ then $$\s(ht)=(\s\cdot
t^{-1})(h)=(\s\cdot t^{-1})(1)=\s(t).$$ Thus $\sigma$ is constant on
the coset $Ht$ and we may regard $S$ as a subspace of the space
$C^0(H\backslash G)$ of continuous functions on $H\backslash G.$

We claim that $S$ separates points of  $H\backslash G.$ Indeed suppose 
that $t,t_1\in G$
are such that $\s(Ht)=\s(Ht_1)$ for all $\s\in S.$ Then $\s(t)=\s(t_1),$
so that 
$$\s(tt^{-1}_1)=(\s\cdot t_1)(t)=(\s\cdot t_1)(t_1)=\s(t_1t_1^{-1})=
\s(1).$$
Thus $tt_1^{-1}\in H$ and $Ht=Ht_1,$ establishing the claim. We 
are going to prove

\begin{lem}\label{l3.1} Assume that either $G$ is finite, or that $U$ is 
simple and orthogonally complemented. Then $S={}^HF$ is the space of 
$H$-invariants\footnote{Previously we wrote $V^G$ etc for $G$-invariants, 
even though the action
was a left action. Here, however, we need to distinguish between 
left and right actions of $G$ on $F.$} under the left action of $G$ 
on $F.$
\end{lem}

If this is true then it follows that $R={}^HW.$ Then clearly $U$ is the
space of $H$-invariants on $V,$ and the main theorem is proved.

There are natural identifications which give containments $S\subset
{}^HF\subset C^0(H\backslash G),$ and we also know that $S$ separates points of
$H\backslash G.$ If $G$ is finite then this already forces $S=C^0(G/H)$ and
there is nothing more to prove. So part (B) of  Theorem \ref{thmm} is
established.

To complete the proof of the lemma we show that $S$ is dense in
$C^0(H\backslash G)$ (with the supremum norm). Since ${}^HF$ is a
direct sum of finite-dimensional subspaces ${}^H\End_{\C}(W_{\l}),$
each of which is closed, we then conclude immediately that $S={}^HF,$
as required.  According to the Stone-Weierstrass theorem, the density
of $S$ in $C^0(H\backslash G)$ will follow as long as we can show that
$S$ is closed under complex conjugation, since we have already seen that
$S$ separates points of $H\backslash G.$  It is this condition which
requires us to assume -- as we now do -- that $U$ is also simple and 
orthogonally complemented.

For $\l\in \L,$ let $\lambda^*$ be the dual of $\l$ in the sense that 
$W_{\l^*}=W_{\l}^*=\Hom_{\C}(W_{\l},\C).$ We restrict the canonical
projection $W_{\l^*}\otimes W_{\l}\to W_1=\C$ to
a pairing $R_{\l^*}\otimes R_{\l}\to \C,$ and prove next 
that this latter pairing is non-degenerate.

For convenience we set $U_{\l}=R_{\l}\otimes V_{\l}.$ Since $U$ is a simple
VOA it follows from  Proposition 4.1 of [DM2] that
for any nonzero $v\in R_{\l}$ we have
$$U=\<u_mv|u\in U,m\in\Z\>.$$ 
Notice that $u_mv\in \oplus_{\nu\in \mu\otimes \l}U_{\nu}$ if $u\in U_{\mu}.$
Thus $\<u_mv|u\in U_{\mu}, \mu\in\L\setminus \{\l^*\}, m\in \Z\>$ is contained in
$\oplus_{\nu\ne 1}U_{\nu}$ and $V^G\subset \<u_mv|u\in U_{\l^*},m\in \Z\>.$
Since $U_{\l^*}$ is a $V^G$-module generated by $R_{\l^*}$ we see from
the associativity of vertex operators that
$\<u_mv|u\in U_{\l^*},m\in \Z\>$ is a $V^G$-module generated
by  $\<u_mv|u\in R_{\l^*},m\in \Z\>.$ The actions of $u_m$ and $G$ 
on $V$ commute for $u\in V^G$ and $m\in \Z.$ It follows 
that there exists $u\in R_{\l^*}$ and $p\in \Z$ such that the projection
of $u_pv$ into $V^G$ is nonzero. 

Consider the composition of $G$-maps
$$W_{\l^*}\otimes W_{\l}\stackrel{\psi_p}{\longrightarrow}Z_p\stackrel{\tau}{\longrightarrow}Z_p^G$$
where we use the notation of (\ref{2.9})-(\ref{2.10}) and $\tau$ is projection
onto $G$-invariants. From the last paragraph it follows that the image of 
$R_{\l^*}\otimes v$ under $\tau\circ \psi_p$ is nonzero. Thus $\tau\circ\psi_p$ generates $\Hom_{G}(W_{\l^*}\otimes W_{\l},\C)$ and its restriction to 
$R_{\l^*}\otimes v$ is nonzero. Since $v$ is an arbitrary non-zero element of $R_{\l},$ and since the same argument applies with $\l$
and $\l^*$ interchanged, this proves that indeed restriction of the canonical pairing
to $R_{\l^*}\otimes R_{\l}$ is nondegenerate.

Next we claim that the restriction of the pairing $W_{\l^*}\otimes W_{\l}\to\C$
to both $R^c_{\l^*}\otimes R_{\l}$ and  $R_{\l^*}\otimes R_{\l}^c$
is zero. The proofs of these two assertions being similar we
only prove the second. Namely, if $R_{\l^*}\otimes R_{\l}^c\to\C$ in 
nonzero, application of Lemma \ref{l2.3} shows that 
that there are $u\in R_{\l^*},$ $v\in  R_{\l}^c$
and $m\in\Z$ such that $u_mv$ has a nonzero projection into $V^G.$
But $u_mv$ lies in $U^c$ and $U^c\cap V^G=0.$ This
contradiction proves the assertion.

Now let us identify $W_{\l^*}$ with $W_{\l}$ via a positive definite,
$G$-invariant, hermitian form on $W_{\l}.$ It follows from what we have
established so far, together with the fact that $R_{\l}$ and
$R_{\l}^c$ are mutually orthogonal, as are $R_{\l^*}$ and
$R_{\l^*}^c,$ that $R_{\l^*}$ is then identified with $R_{\l}.$ In terms
of coordinate functions, this means that if we embed $R_{\l}$ into
$\End\,W_{\l}$ as the first $k$-rows of the first column, i.e., as the 
coordinate functions $a_{i1},$ $1\leq i\leq k,$ then
$R_{\l^*}$ corresponds to the coordinate functions $\bar a_{i1},$
$1\leq i\leq k.$

It follows that for all $\l,$ $\bar R_{\l}\subset S,$ and hence 
$S$ is indeed closed under complex conjugation. This completes 
the proof of Lemma \ref{l3.1} and with it that of Theorem \ref{thmm}. 

\section{Further comments}
\setcounter{equation}{0}

First we wish to point out that the statement of Theorem 3.1 of [DLM]
inadvertently omitted the adjective ``simple'' so that as stated, 
Theorem 3.1 of (loc. cit.) is incorrect. Of course, the correct statement
is contained as a special case of the main theorem of the present paper.
Note that in the situation of (loc.cit.), where $G$ is abelian, the condition
of orthogonal complementation is redundant.
As we pointed out in the introduction, Galois correspondences of the type
discussed in this paper must necessarily only involve simple sub VOAs,
although there are generally many sub VOAs $U$ with
$V^G\subset U\subset V$ and $U$ not simple, at least if $G$
is not finite.

Indeed, suppose that $\L$ is as before, and that $\Theta\subset \L$ is a 
{\em closed} subset in the sense that if $\alpha,\b\in \Theta,$ then
$\delta\in \Theta$ whenever $W_{\delta}$ is a constituent of 
$W_{\l}\otimes W_{\b}.$ Suppose that $1\in \Theta.$ Then 
\begin{equation}\label{4.1}
U=\bigoplus_{\l\in\Theta}W_{\l}\otimes V_{\l}
\end{equation}
is a sub VOA of $V$ which contains $V^G.$ What we have proved implies that 
$U$ is a simple sub VOA if, and only if, $\Theta$ is closed with
respect to conjugation, i.e., if $\l\in \Theta$ then also $\l^*\in \Theta.$
Finally, note that  (\ref{4.1}) represents a typical $G$-invariant sub VOA of $V$ which contains $V^G,$ so that there is a bijection between $G$-invariant
sub VOAs containing $V^G$ and closed subsets of $\L$ containing 1.

We consider a different kind of example. Let $V_L$ be the simple
vertex operator algebra associated to the lattice $L=\Z\alpha$ where
$\<\alpha,\alpha\>=2n$ for positive integer $n$ (cf. [B] and
[FLM]). Let $L^{\circ}=\frac{1}{2n}\Z\a$ be the dual lattice of $L.$
Then the compact Lie group $G=R\alpha/L^{\circ},$ which is a circle,
acts on $V_L$ in the following way (cf. [DM1]): $$ \beta\cdot
(u\otimes e^{\gamma})=e^{2\pi i\<\beta,\gamma\>} u\otimes e^{\gamma}$$
where $\beta\in G,$ $u\in \C[\alpha(-1),\alpha(-2),\cdots ]$ and
$\gamma\in L.$ It is easy to see that
$$V_L^G=\C[\alpha(-1),\alpha(-2),\cdots ].$$

Let $U=\oplus_{m\geq 0}\C[\alpha(-1),\alpha(-2),\cdots ]\otimes e^{m\alpha}.$
Then $U$ is a sub VOA of $V_L$ which contains $V_L^G.$ Note that
$U$ is not simple as $\oplus_{m\geq 1}\C[\alpha(-1),\alpha(-2),\cdots ]\otimes e^{m\alpha}$ is a proper ideal of $U.$ 
Clearly if
$g\in G$ is the identity on $U$ then $g$ must be the identity on the whole
space $V_L.$ Thus $U\ne V^H_L$ for any subgroup $H$ of $G.$ This kind of 
sub VOA was also considered in [D] in order to construct examples of vertex
operator algebras with certain properties.

\end{document}